\documentclass[10pt]{article}
\usepackage[final]{epsfig}
\usepackage{graphics}
\usepackage{amsmath}
\usepackage{amsfonts}
\usepackage{latexsym}
\usepackage{amssymb}
\usepackage{graphicx}
\usepackage{url}

\newtheorem{theorem}{Theorem}

\newcommand{\eps}{{\varepsilon}}
\newcommand{\g}{{\gamma}}
\newcommand{\G}{{\Gamma}}

\newcommand{\RP}{{\mathbf {RP}}}

\begin{document}

\title{Osculating curves: around the Tait-Kneser Theorem}

\author{E. Ghys \and  S. Tabachnikov \and V. Timorin}

\date{}
\maketitle

\section{Tait and Kneser} \label{TK}

The notion of osculating circle (or circle of curvature) of a smooth plane curve is familiar to every student of calculus and elementary differential geometry: this is the circle that approximates the curve at a point better than all other circles. 

One may say that the osculating circle passes through three infinitesimally close points on the curve. More specifically, pick three points on the curve and draw a circle through these points. As the points tend to each other, there is a limiting position of the circle: this is the osculating circle. Its radius is the radius of curvature of the curve, and the reciprocal of the radius is the curvature of the curve. 

If both the curve and the osculating circle  are represented locally as graphs of smooth functions then not only the values  of these functions but also their first and second derivatives coincide at the point  of contact.

\begin{figure}[hbtp]
\centering
\includegraphics[width=1.2in]{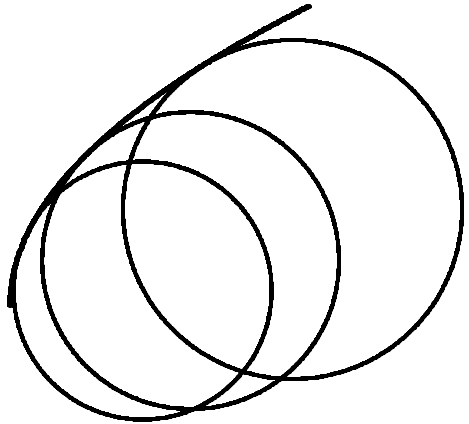}
\caption{Osculating circles?}
\label{wrong}
\end{figure}

Ask your mathematical friend to sketch an arc of a curve and a few osculating circles. Chances are, you will see something like Figure \ref{wrong}.

This is wrong! The following theorem was discovered by Peter Guthrie Tait in the end of the 19th century \cite{Tai} and rediscovered by Adolf Kneser  early in the 20th century \cite{Kne}.

\begin{theorem} \label{Tait}
The osculating circles of an arc with monotonic positive curvature are pairwise disjoint and nested.
\end{theorem}

 Tait's paper is so short that we quote it almost verbatim (omitting some old-fashioned terms):
\begin{quote}
\small {When the curvature of a plane curve continuously increases or diminishes (as in the case with logarithmic spiral for instance) no two of the circles of curvature can intersect each other. 

This curious remark occurred to me some time ago in connection with an accidental feature of a totally different question...

The proof is excessively simple. For if $A, B$, be any two points of the evolute, the chord $AB$ is the distance between the centers of two of the circles, and is necessarily less than the arc $AB$, the difference of their radii...

When the curve has points of maximum or minimum curvature, there are corresponding $\dots$ cusps on the evolute; and pairs of circles of curvature whose centers lie on opposite sides of the cusp, $C$, may intersect: -- for the chord $AB$ may now exceed the difference between $CA$ and $CB$.}
\end{quote}

\begin{figure}[hbtp]
\centering
\includegraphics[width=2.8in]{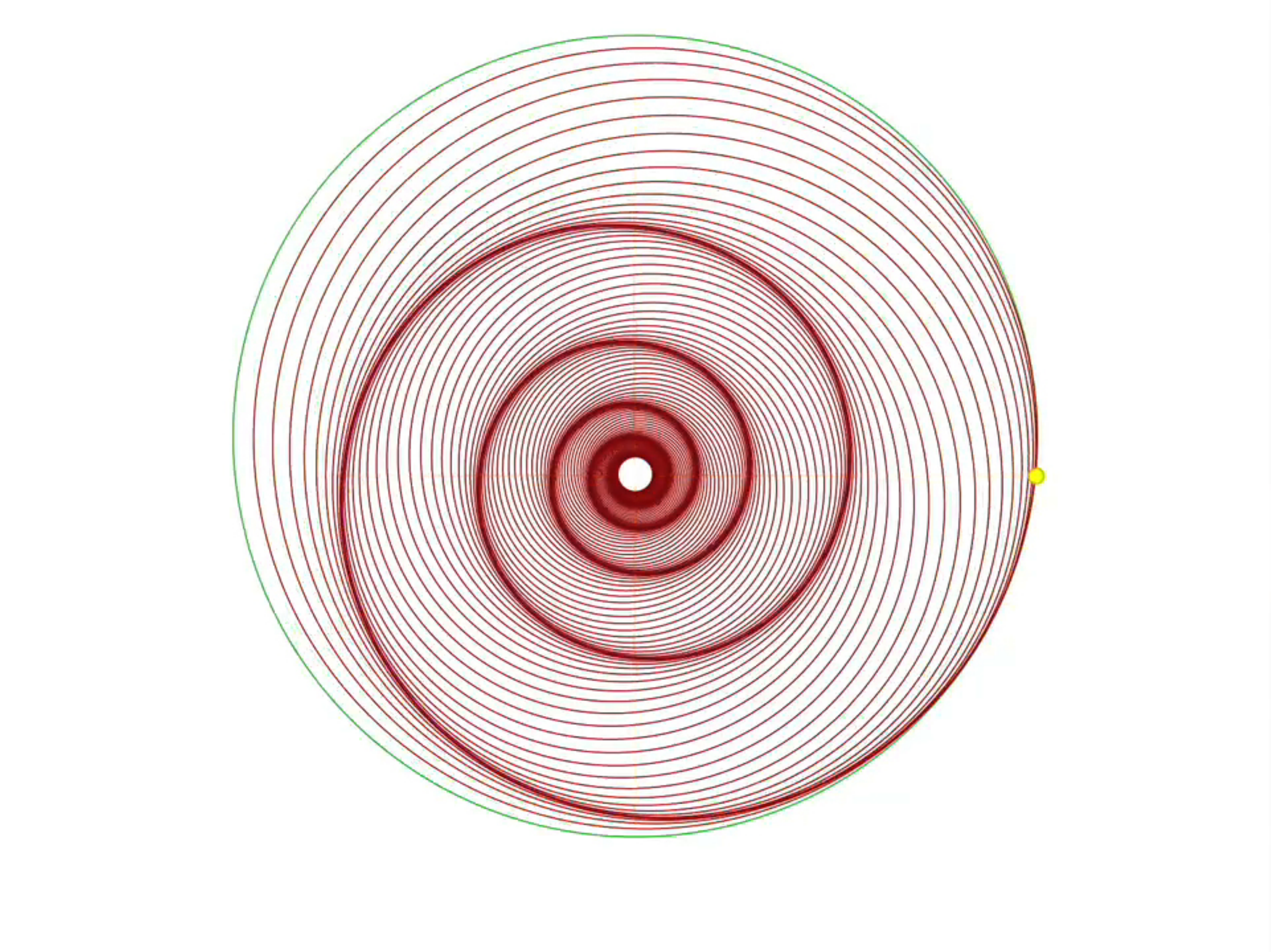}
\caption{Osculating circles of a spiral.}
\label{spiral}
\end{figure}

See Figure \ref{spiral} for a family of osculating circles of a spiral.\footnote{Curiously, the current English Wikipedia article on osculating circles contains three illustrations, and none of them depicts the typical situation: the curve goes from one side of the osculating circle to the other. The French Wikipedia article fares better in this respect; the reader may enjoy researching other languages.}

\section{Evolutes and involutes} \label{evinv}

Perhaps a hundred years ago Tait's argument was self-evident and did not require further explanation. Alas, the situation is different today, and this section is an elaboration of his proof. The reader is encouraged to consult her favorite book on elementary differential geometry for the basic facts that we recall below.

\begin{figure}[hbtp]
\centering
\includegraphics[height=2in]{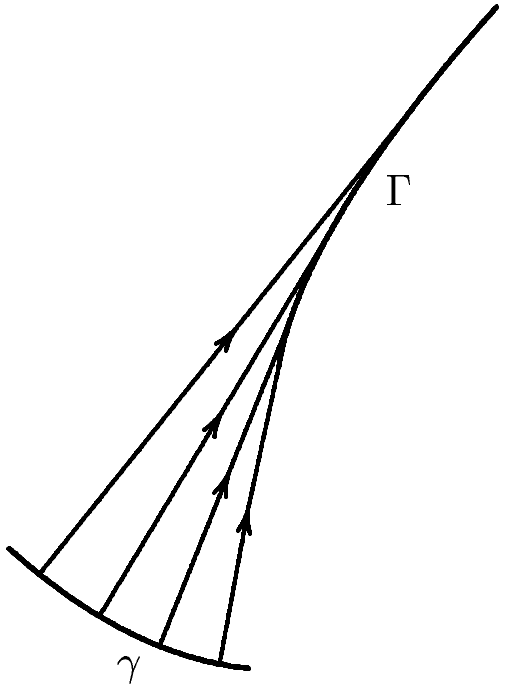}
\quad\quad
\includegraphics[height=2in]{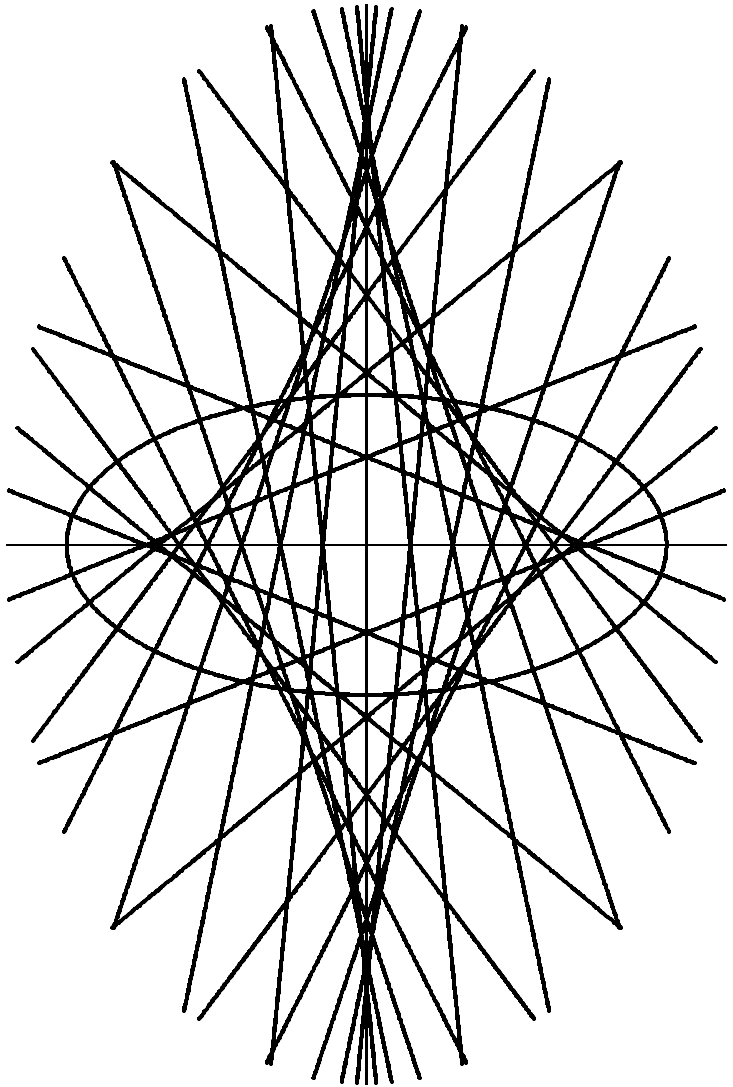}
\caption{$\G$ is the evolute of $\g$. The evolute of an ellipse.}
\label{evolute}
\end{figure}

The locus of centers of osculating circles is called the {\it evolute} of a curve. The evolute is also the envelope of the family of normal lines to the curve. See Figures \ref{evolute}.

The evolute typically has cusp singularities, clearly seen in Figure \ref{evolute}. For generic curves, these are the centers of the stationary osculating circles, the osculating circles at the {\it vertices} of the curve, that is, the points where the curvature has a local minimum or a local maximum.

Consider the left Figure \ref{evolute} again. The curve $\g$ is called an {\it involute} of the curve $\G$: an involute  is orthogonal to the tangent lines of a curve. The involute $\g$ is described by the free end of a non-stretchable string whose other end is fixed on   $\G$ and which is wrapped around it (for this reason, involutes are also called evolvents). That this {\it string construction} indeed does the job is obvious: the radial component of the velocity of the free end point would stretch the string.

A consequence of the string construction is that the length of an arc of the evolute $\G$ equals the difference of its 
tangent segments to the involute $\g$, that is, the increment of the radii of curvature of $\g$. This is true as long as the curvature of $\g$ is monotonic and $\G$ is free of cusps.

Another curious consequence is that the evolute of a closed curve has total length zero. The length is algebraic: its sign changes each time that one passes a cusp. We leave it to the reader to prove this zero length property (necessary and sufficient for the string construction to yield a closed curve). 

\begin{figure}[hbtp]
\centering
\includegraphics[width=1.5in]{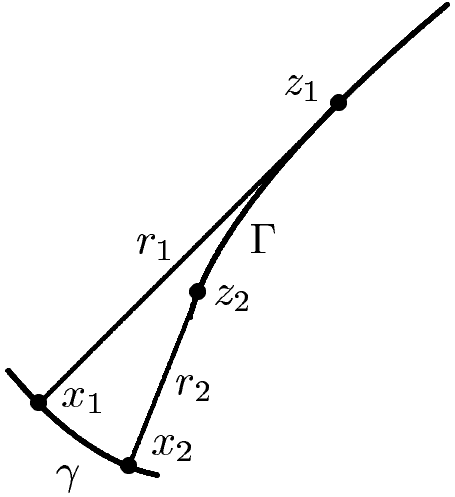}
\caption {Tait's proof.}
\label{PfTK}
\end{figure}

Tait's argument is  straightforward now, see Figure \ref{PfTK}. Let $r_1$ and $r_2$ be the radii of osculating circles at points $x_1$ and $x_2$, and $z_1$ and $z_2$ be their centers. Then the length of the arc $z_1 z_2$  equals $r_1-r_2$, hence $|z_1 z_2| < r_1-r_2$. Therefore the circle with center $z_1$ and radius $r_1$ contains the circle with center $z_2$ and radius  $r_2$.

\section{A paradoxical foliation} \label{fol}

Let us take a look  at Figure \ref{spiral} again. We see an annulus bounded by the smallest and the largest of the osculating circles of a curve $\g$ with monotonic curvature. This annulus is foliated by the osculating circles of $\g$ , and the curve ``snakes" between these circles, always remaining tangent to them. How could this be possible? 

Isn't this similar to having a non-constant function with everywhere zero derivative?
Indeed, if the foliation consists of horizontal lines and the curve is the graph of a differentiable function $f(x)$, then $f'(x)=0$ for all $x$, and $f$ is constant. But then the curve is contained within one leaf. 

The resolution of this ``paradox" is that this foliation is not differentiable and we cannot locally map the family of osculating circles to the family of parallel lines by a smooth map. A foliation is determined by a function whose level curves are the leaves; a foliation is differentiable if this function can be chosen differentiable. A foliation may have  leaves as good as one wishes (smooth, analytic, algebraic) and still fail to be differentiable.

\begin{theorem} \label{foliation}
If a differentiable  function in the annulus is constant on each osculating circle then this is a constant  function.
\end{theorem}

For example, the radius of a circle is a function constant on the leaves. As a function in the annulus, it is not differentiable.

To prove the theorem, let $F$ be a differentiable function constant on the leaves. Then  $dF$ is a differential 1-form whose restriction to each circle is zero. The curve $\g$ is tangent to one of these circles at each point. Hence $dF$ is zero on $\g$ as well. Therefore  $F$ is constant on $\g$. But $\g$ intersects all the leaves, so $F$ is constant in the annulus. 

Thus a perfectly smooth (analytic, algebraic) curve provides an example of a non-differentiable foliation by its osculating circles. 

\section{Taylor polynomials} \label{Taylor}

In this section we present a version of Tait-Kneser theorem for Taylor polynomials. It is hard to believe that this result was not known for a long time, but we did not see it in the literature.

Let $f(x)$ be a smooth function of real variable.  The Taylor polynomial $T_t(x)$ of
degree $n$ approximates $f$ up to the $n$-th derivative:
$$
T_t (x) =\sum_{i=0}^n \frac{f^{(i)}(t)}{i!} (x-t)^i.
$$
Assume that $n$ is even and that $f^{(n+1)} (x)\neq0$ on some interval
$I$.

\begin{theorem} \label{TKpol}
For any distinct $a,b \in I$, the graphs of the Taylor
polynomials $T_a$ and $T_b$ are disjoint over the whole real line.
\end{theorem}

To prove this, assume that $f^{(n+1)} (x)>0$ on $I$ and that $a<b$. 
One has:
$$
\frac{\partial T_t}{\partial t} (x)= \sum_{i=0}^n \frac{f^{(i+1)}(t)}{i!} (x-t)^i - \sum_{i=0}^n \frac{f^{(i)}(t)}{(i-1)!} (x-t)^{i-1} = \frac{f^{(n+1)}(t)}{n!} (x-t)^n,
$$
and hence $(\partial T_t/\partial t) (x)> 0$ (except for $x=t$). It follows that $T_t(x)$ increases, as a function of $t$, therefore  $T_a(x)<T_b(x)$ for all $x$.

The same argument proves the following variant of Theorem
\ref{TKpol}. Let $n$ be odd, and assume that  $f^{(n+1)} (x)\neq0$ on
an interval $I$. 

\begin{theorem} \label{TKpol1}
For any distinct $a,b \in I,\ a<b$, the graphs of the Taylor
polynomials $T_a$ and $T_b$ are disjoint over the interval $[b,\infty)$.
\end{theorem}

\begin{figure}[hbtp]
\centering
\includegraphics[height=1.1in]{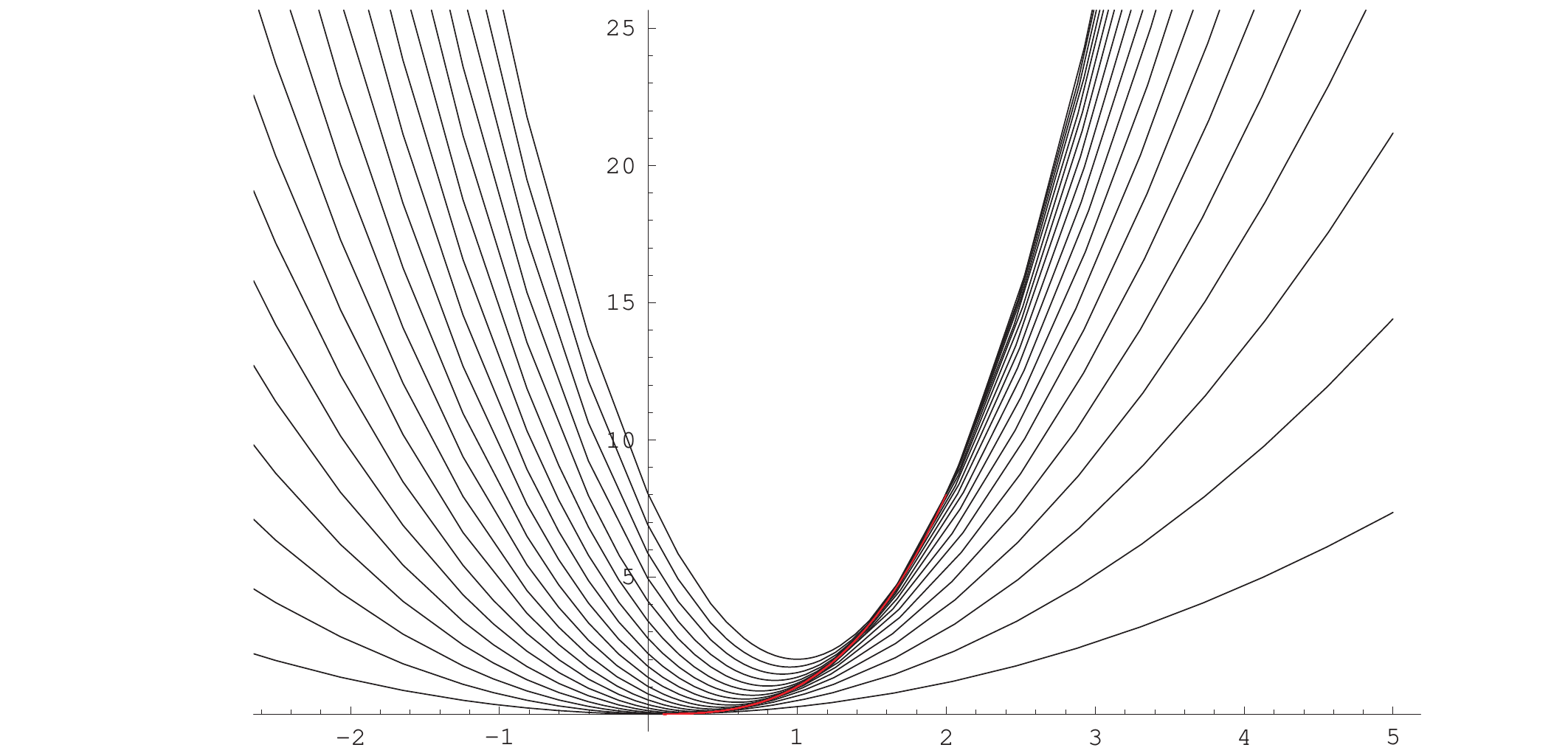}
\ 
\includegraphics[height=1.1in]{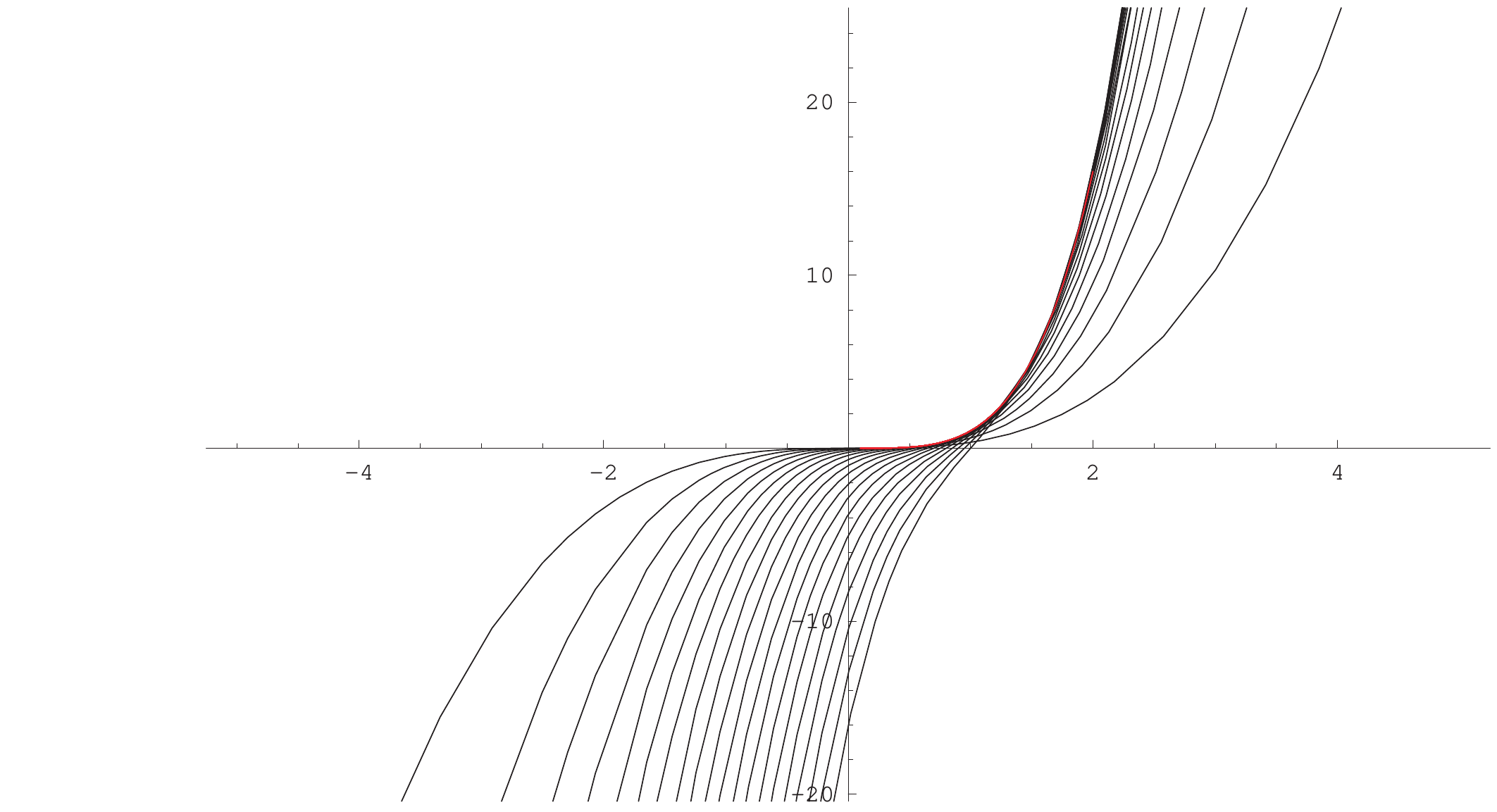}
\caption{Quadratic Taylor polynomials of the function $f(x)=x^3$ and cubic Taylor polynomials of the function $f(x)=x^4$}
\label{2and3}
\end{figure}

Theorems \ref{TKpol} and \ref{TKpol1} are illustrated in Figure \ref{2and3}.

The same proof establishes more: not only the function $T_b(x) -T_a(x)$ is positive, but it is also convex. Furthermore, all its derivatives of even orders are positive. Certain analogs of this remark apply to the variations on the Tait-Kneser theorem presented in the next section, but we shall not dwell on this intriguing subject here.

\section{Variations} \label{var}

The Tait-Kneser theorem can be extended from  circles to other classes of curves. Let us consider a very general situation when a $d$-parameter family of plane curves is given; these curves will be used to approximate a test smooth curve at a point. For example, a conic depends on five parameters, so $d=5$ for the family of conics.

Given a smooth curve $\g$ and point $x\in\g$, the osculating curve from our family is the curve that has tangency with $\g$ at point $x$ of order $d-1$; in other words, it is the curve from the family that passes through $d$ infinitesimally close points on $\g$. The curve hyperosculates if the order of tangency is greater, that is, the curve passes through $d+1$ infinitesimally close points on $\g$.   

For example, one has the 1-parameter family of osculating conics of a plane curve $\g$ parameterized by the point $x\in\g$. A point $x$ is called {\it sextactic} if the osculating conic  hyperosculates at this point. In general, a point of $\g$ is called {\it extactic} if the osculating curve  hyperosculates at this point.

We shall now describe a number of Tait-Kneser-like theorems. Our discussion is informal; the reader interested in more details is refereed to \cite{Gh2,TT}. Let us consider the case of osculating conics.

\begin{figure}[hbtp]
\centering
\includegraphics[width=2in]{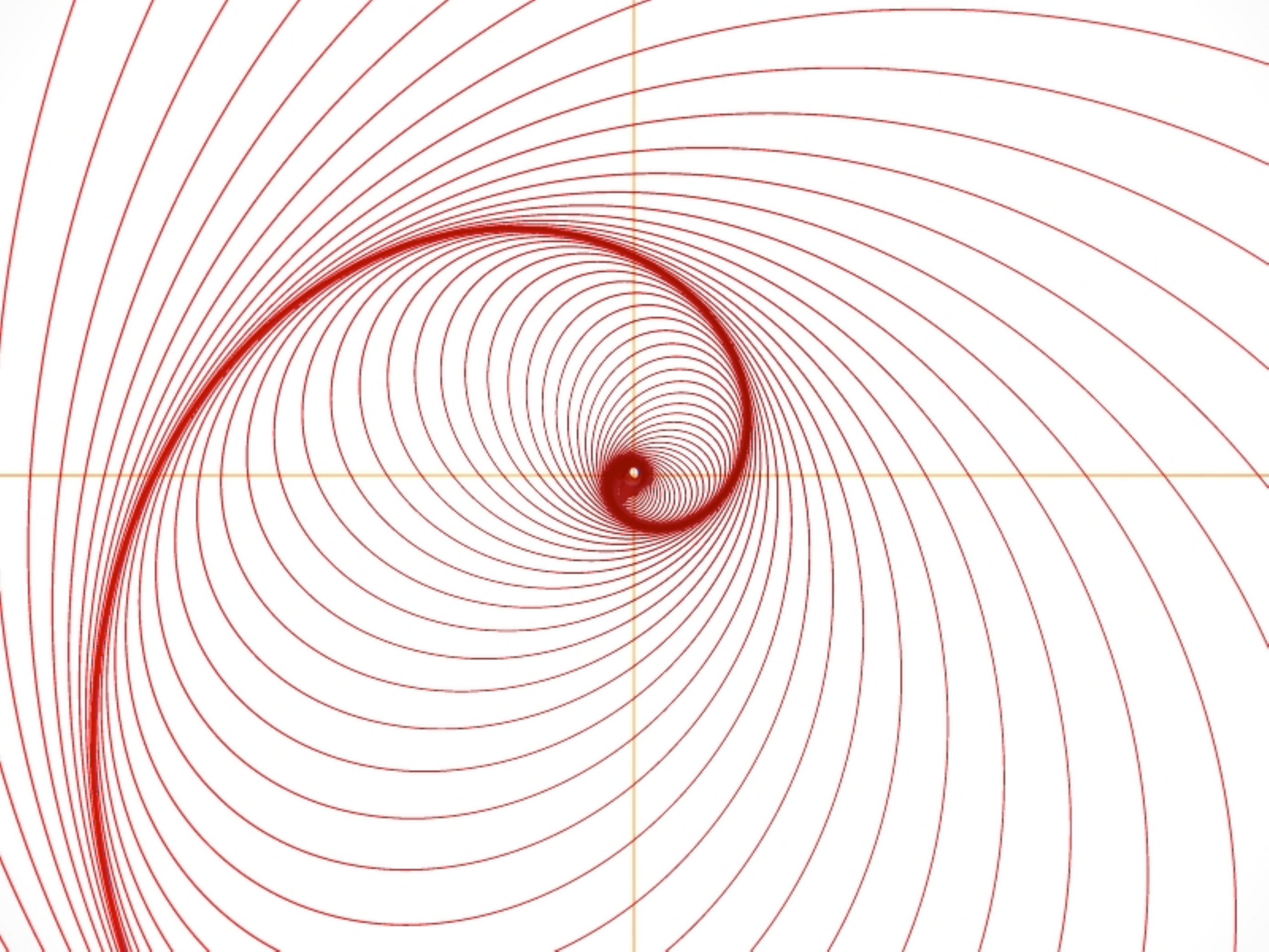}
\caption {Osculating conics of a spiral.}
\label{oscconics}
\end{figure}

\begin{theorem} \label{conics}
The osculating conics of a curve, free from sextactic points, are pairwise disjoint and nested (see Figure \ref{oscconics}).
\end{theorem}

This theorem is better understood in the projective plane where all non-degenerate conics are equivalent, and there is no difference between ellipses, parabolas and hyperbolas. In particular, a non-degenerate conic divides the projective plane into two domains, the inner one which is a disc, and the outer one which is the M\"obius band.

Here is a sketch of a proof.\footnote{A similar argument applies to osculating circles as well.} Give the curve a parameterization, $\g(x)$, and let $C_x$ be the osculating conic at point $x$. Let $F_x=0$ be a quadratic equation of the conic $C_x$.

It suffices to establish the claim for sufficiently close osculating conics, so consider  infinitesimally close ones. The intersection of the conics $C_x$ and $C_{x+\eps}$ (for infinitesimal $\eps$) is given by the system of equations
$$
F_x=0,\quad \frac{\partial F_x}{\partial x}=0.
$$
Both equations are quadratic so, by the Bezout theorem, the number of solutions is at most 4 (it is not infinite because $x$ is not a sextactic point). But the conics $C_x$ and $C_{x+\eps}$ already have an intersection of multiplicity 4 at point $x$: each is determined by 5 ``consecutive" points on the curve $\g$, and they share 4 of these points. Therefore they have no other intersections, as needed. 	

Another generalization, proved similarly, concerns diffeomorphisms of the real projective line $\RP^1$. At every point, a diffeomorphism $f: \RP^1\to \RP^1$ can be approximated, up to the second derivative, by a fractional-linear (M\"obius) transformation 
$$
x\mapsto \frac{ax+b}{cx+d}.
$$
It is natural to call this the osculating M\"obius transformation of $f$. Hyperosculation occurs when the approximation is finer, up to the third derivative; this happens when the Schwarzian derivative of $f$ vanishes:
$$
S(f)(x)=\frac{f'''(x)}{f'(x)}-\frac{3}{2} \left(\frac{f''(x)}{f'(x)}\right)^2 =0
$$
(see \cite{OT,OT1} concerning the Schwarzian derivative).

\begin{theorem} \label{Sch}
Let $f:[a,b]\to \RP^1$ be a local diffeomorphism whose Schwarzian derivative does not vanish. Then the graphs of the osculating M\"obius transformation are pairwise disjoint.
\end{theorem}

Of course, these graphs are hyperbolas with the vertical and horizontal asymptotes.

Can one generalize to  algebraic curves of higher degree? The space of algebraic curves of degree $d$ has dimension $n(d)=d(d+3)/2$.
The osculating algebraic curve of degree $d$  passes through $n(d)$ infinitesimally close points of a smooth curve $\g$. Two infinitesimally close osculating curves of degree $d$ at point $x\in\g$ have there an intersection of multiplicity $n(d)-1$, whereas two curves of degree $d$ may have up to $d^2$ intersections altogether. For $d\geq 3$, one has $d^2>d(d+3)/2 -1$, so one cannot exclude intersections of osculating algebraic curves of degree $d$.

\begin{figure}[hbtp]
\centering
\includegraphics[width=3in]{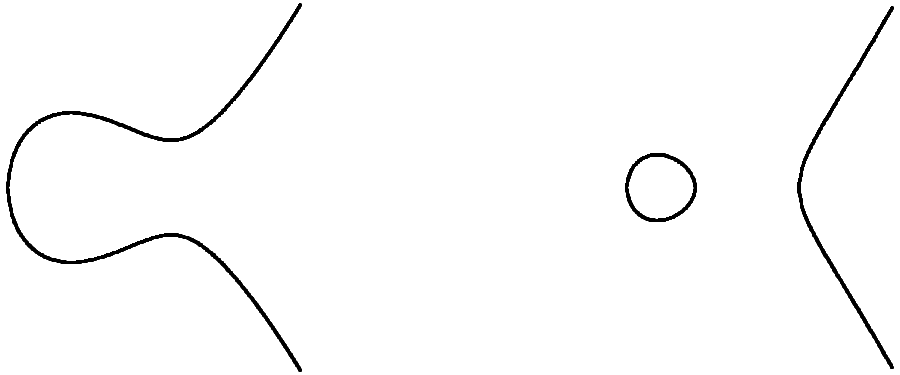}
\caption {Two types of cubic curves.}
\label{cubiccurve}
\end{figure}

However, one can remedy the situation  for cubic curves. A cubic curve looks like shown in Figure \ref{cubiccurve}: it may have one or two components, and in the latter case one of them is compact. The compact component is called the oval of a cubic curve. Two ovals intersect in an even number of points, hence one can reduce the number $9=3^2$ to 8 if one considers ovals of cubic curves as osculating curves. This yields

\begin{theorem} \label{cube}
Given a plane curve, osculated by ovals of cubic curves and free from extactic points, the osculating ovals are disjoint and pairwise nested.
\end{theorem}

\begin{figure}[hbtp]
\centering
\includegraphics[width=2.5in]{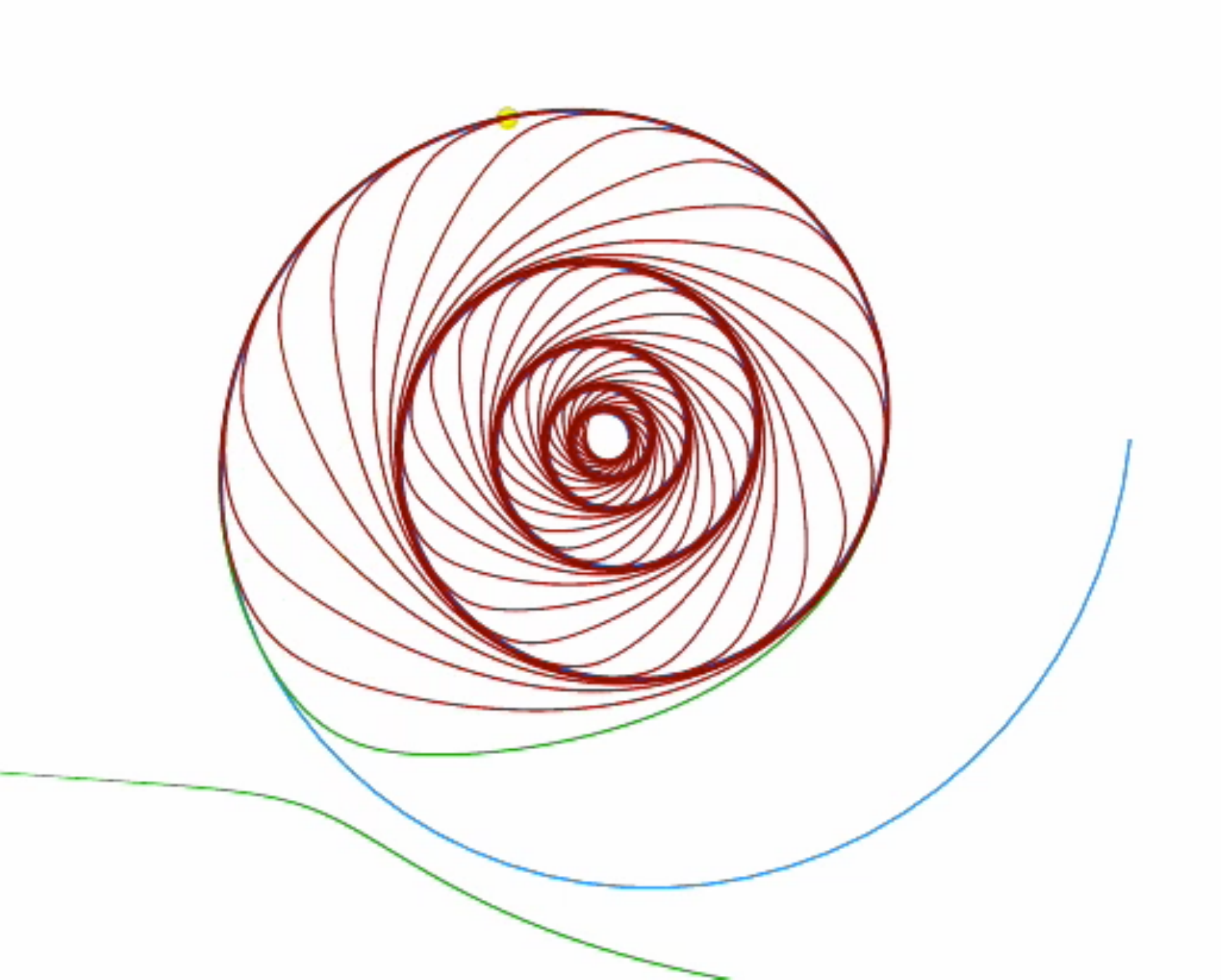}
\caption {A spiral osculated by ovals of cubic curves.}
\label{cubics}
\end{figure}

See Figure \ref{cubics} for an illustration.

\section{4-vertex theorem and beyond} \label{four}

This story would be incomplete without mentioning a close relation of various versions of the Tait-Kneser theorem and numerous results on the least number of extactic points. The first such result is the 4-vertex theorem discovered by S. Mukhopadhyaya in 1909 \cite{Mu}: {\it a plane oval \footnote{Closed smooth strictly convex curve} has at least four vertices}. In the same paper, Mukhopadhyaya proved the 6-vertex theorem: {\it a plane oval has at least six sextactic points}. Note that these numbers, 4 and 6, are one greater than the dimensions of the respective spaces of osculating curves, circles and conics.

A similar theorem holds for M\"obius transformations approximating diffeomorphisms of the projective line: {\it for every diffeomorphism of $\RP^1$, the Schwarzian derivative vanishes at least four times} \cite{Gh1}. 

And what about approximating by cubic curves? Although not true for arbitrary curves, the following result holds: {\it a plane oval, sufficiently close to an oval of a cubic curve, has at least 10 extactic points} \cite{Ar}. Once again, $10=9+1$ where 9 is the dimension of the space of cubic curves. We refer to \cite{OT} for information about the 4-vertex theorem and its relatives.

By the way, the reader may wonder whether there is a ``vertex" counterpart to Theorem \ref{TKpol}. Here is a candidate: {\it if $f(x)$ is a smooth function of real variable, flat at infinity (for example, coinciding with $\exp({-x^2})$ outside of some interval), then, for each $n$, the equation $f^{(n)}(x)=0$ has at least $n$ solutions.} The proof easily follows from the Rolle theorem.

One cannot help wondering about the meaning of this relation between two sets of theorems. Is there a general underlining principle in action here?

\bigskip
{\bf Acknowledgments}. We are grateful to Jos Leys for producing images used in this article. 
 S. T. was partially supported by the Simons Foundation grant No 209361 and by the NSF grant DMS-1105442. V. T. was 
partially supported by the Deligne fellowship, the Simons-IUM fellowship, RFBR grants 10-01-00739-a,
11-01-00654-a, MESRF grant MK-2790.2011.1, and AG Laboratory NRU-HSE, MESRF grant ag. 11 11.G34.31.0023.

\end{document}